\documentclass[onefignum,onetabnum]{siamart190516}
\pdfoutput=1

\usepackage{amsfonts}
\usepackage{graphicx}
\graphicspath{{Figures/}}
\usepackage{epstopdf}
\usepackage{tabularx}
\usepackage{multirow}
\usepackage{booktabs}
\usepackage{mathtools}
\usepackage{algpseudocode}
\usepackage[english]{babel}
\usepackage{enumitem} 
\usepackage{amsopn}

\ifpdf
  \DeclareGraphicsExtensions{.eps,.pdf,.png,.jpg}
\else
  \DeclareGraphicsExtensions{.eps}
\fi

\newsiamremark{remark}{Remark}
\newsiamremark{hypothesis}{Hypothesis}
\crefname{hypothesis}{Hypothesis}{Hypotheses}
\newsiamthm{claim}{Claim}

\headers{Semi-Structured Algebraic Multigrid}{V. A. P. Magri, R. D. Falgout, and U. M. Yang}

\title{A New Semi-Structured Algebraic Multigrid Method
       \thanks{Submitted to the editors on July 15, 2021.
       {\protect\\ \hspace*{10.37pt} This work was performed under the auspices of the
         U.S. Department of Energy by Lawrence Livermore National Laboratory under Contract DE-AC52-07NA27344. LLNL-JRNL-834288-DRAFT.}}}

\author{Victor A. P. Magri\thanks{Lawrence Livermore National Laboratory
  (\email{paludettomag1@llnl.gov}).}
\and Robert D. Falgout\thanks{Lawrence Livermore National Laboratory
  (\email{falgout2@llnl.gov}).}
\and Ulrike M. Yang\thanks{Lawrence Livermore National Laboratory
  (\email{yang11@llnl.gov}).}
}

\ifpdf
\hypersetup{
  pdftitle={A New Semi-Structured Algebraic Multigrid Method},
  pdfauthor={V. A. P. Magri, R. D. Falgout, and U. M. Yang}
}
\fi
\providecommand{\mat}[1]{#1}
\providecommand{\matT}[1]{\mat{#1}^{T}}
\providecommand{\matI}[1]{\mat{#1}^{-1}}

\renewcommand{\vec}[1]{\mathbf{#1}}

\def\hypre{\textsl{hypre}}

\algnewcommand\algorithmicbloopb{\textbf{BoxLoopBegin}}
\algdef{S}[FOR]{BoxLoopBegin}[2]{\algorithmicbloopb\ (#1, #2)}
\algnewcommand\algorithmicbloope{\textbf{BoxLoopEnd}}
\algdef{E}[FOR]{BoxLoopEnd}{\algorithmicbloope}
\algnewcommand\algorithmicbloopbb{\textbf{BoxLoopBegin2}}
\algdef{S}[FOR]{BoxLoopBeginB}[3]{\algorithmicbloopbb\ (#1, #2, #3)}
\algnewcommand\algorithmicbloopbe{\textbf{BoxLoopEnd2}}
\algdef{E}[FOR]{BoxLoopEndB}{\algorithmicbloopbe}

\DeclareMathOperator*{\argminD}{arg\,min_d}
\DeclareMathOperator*{\vol}{volume}

\providecommand{\ParCSRMatrix}{\texttt{ParCSRMatrix}}
\providecommand{\IJ}{\texttt{IJ}}
\providecommand{\ParCSR}{\texttt{ParCSR}}
\providecommand{\Struct}{\texttt{Struct}}
\providecommand{\SStruct}{\texttt{SStruct}}
\providecommand{\SStructSplit}{\text{Split}}
\providecommand{\Split}{\text{Split}}
\providecommand{\SStructGrid}{\texttt{SStructGrid}}

\providecommand{\SStructMatrix}{\texttt{SStructMatrix}}

\providecommand{\SStructSplitOpt}{\text{Split}}
\providecommand{\BoomerAMGOpt}{\text{BoomerAMG}}
\providecommand{\PFMGOpt}{\text{PFMG}}
\providecommand{\SSAMGOpt}{\text{SSAMG-opt}}
\providecommand{\SSAMGBase}{\text{SSAMG-base}}
\providecommand{\SSAMGSkip}{\text{SSAMG-skip}}
\providecommand{\SSAMGHybrid}{\text{SSAMG-hybrid}}

\begin{document}

\maketitle

\begin{abstract}
  Multigrid methods are well suited to large massively parallel computer architectures
  because they are mathematically optimal and display excellent parallelization
  properties. Since current architecture trends are favoring regular compute patterns to
  achieve high performance, the ability to express structure has become much more
  important. The \hypre{} software library provides high-performance multigrid
  preconditioners and solvers through conceptual interfaces, including a semi-structured
  interface that describes matrices primarily in terms of stencils and logically
  structured grids. This paper presents a new semi-structured algebraic multigrid (SSAMG)
  method built on this interface. The numerical convergence and performance of a CPU
  implementation of this method are evaluated for a set of semi-structured problems. SSAMG
  achieves significantly better setup times than \hypre{}'s unstructured AMG solvers and
  comparable convergence. In addition, the new method is capable of solving more complex
  problems than \hypre{}'s structured solvers.
\end{abstract}



\section{Introduction}
\label{sec:intro}

The solution of partial differential equations (PDEs) often involves solving linear
systems of equations
\begin{equation}
  \mat{A} \vec{x} = \vec{b},
\label{eq:linsys}
\end{equation}
where $\mat{A} \in \mathbb{R}^{N \times N}$ is a sparse matrix; $\vec{b} \in
\mathbb{R}^{N}$ is the right-hand side vector, and $\vec{x} \in \mathbb{R}^{N}$ is the
solution vector. In modern simulations of physical problems, the number of unknowns $N$
can be huge, e.g., on the order of a few billion. Thus, fast solution methods must be used
for Equation \eqref{eq:linsys}.

Multigrid methods acting as preconditioners to Krylov-based iterative solvers are among
the most common choices for fast linear solvers. In these methods, a multilevel hierarchy
of decreasingly smaller linear problems is used to target the reduction of error
components with distinct frequencies and solve \eqref{eq:linsys} with $O(N)$ computations
in a scalable fashion.  There are two basic types of multigrid methods
\cite{BrHeMc00}. Geometric multigrid employs rediscretization on coarse grids, which needs
to be defined explicitly by the user. A less invasive and less problem-dependent approach
is algebraic multigrid (AMG) \cite{St01}, which uses information coming from the assembled
fine level matrix $\mat{A}$ to compute a multilevel hierarchy. The \hypre{} software
library \cite{Hypre,FaYa02} provides high-performance preconditioners and solvers for the
solution of large sparse linear systems on massively parallel computers with a focus on
AMG methods. It features three different interfaces, a structured, a semi-structured, and
a linear-algebraic interface. Its most used AMG method, BoomerAMG \cite{HeYa2002}, is a
fully unstructured method, built on compressed sparse row matrices (CSR). The lack of
structure presents serious challenges to achieve high performance on GPU
architectures. The most efficient solver in \hypre{} is PFMG \cite{AsFa96}, which is
available through the structured interface. It is well suited for implementation on
accelerators, since its data structure is built on grids and stencils, and achieves
significantly better performance than BoomerAMG when solving the same problems
\cite{BaFKY12,FaLSWY21}; however, it is applicable to only a subset of the problems that
BoomerAMG can solve. This work presents a new semi-structured algebraic multigrid (SSAMG)
preconditioner, built on the semi-structured interface, consisting of mostly structured
parts and a small unstructured component. It has the potential to achieve similar
performance as PFMG with the ability to solve more complex problems.

There have been other efforts to develop semi-structured multigrid methods. For example,
multigrid solvers for hierarchical hybrid grids (HHG) have shown to be highly efficient
\cite{BeHu04,BeWHR07,GmRu14,GmRuSt15,KhRu20}. These grids are created by regularly
refining an initial, potentially unstructured grid. Geometric multigrid methods for
semi-structured triangular grids that use a similar approach have also been proposed
\cite{RoGaLi12}. More recently, the HHG approach has been generalized to a semi-structured
multigrid method \cite{MaBeOhTu21}. Regarding applications, there are many examples
employing semi-structured meshes which can benefit from new semi-structured algorithms,
e.g., petroleum reservoir simulation \cite{GaPeWh19}, marine ice sheets modeling
\cite{StDa13}, next-generation weather and climate models \cite{AdFoHa19}, and solid
mechanics simulators \cite{RuAgZhAl21}, to name a few. In addition, software frameworks
that support the development of block-structured AMR applications such as AMReX
\cite{ZhAlBe19,ZhMyGo21} and SAMRAI \cite{HoKo02} can benefit from the development of
solvers for semi-structured problems.

This paper is organized as follows. Section \ref{sec:sstruct} reviews the semi-structured
conceptual interface of \hypre{}, which enables the description of matrices and vectors
that incorporate information about the problem's structure. Section \ref{sec:ssamg}
describes the new semi-structured algorithm in detail. In section \ref{sec:results}, we
evaluate SSAMG's performance and robustness for a set of test cases featuring distinct
characteristics and make comparisons to other solver options available in
\hypre{}. Finally, in section \ref{sec:conclusions}, we list conclusions and future work.

\section{Semi-structured interface in \hypre{}}
\label{sec:sstruct}

The \hypre{} library provides three conceptual interfaces by which the user can define and
solve a linear system of equations: a structured (\Struct{}), a semi-structured
(\SStruct{}) and a linear algebraic (\IJ{}) interface. They range from highly specialized
descriptions using structured grids and stencils in the case of \Struct{} to the most
generic case where sparse matrices are stored in a parallel compressed row storage format
(\ParCSR{}) \cite{FaJoYa06a, FaJoYa06b}. In this paper, we focus on the \SStruct{}
interface~\cite{FaJoYa06a,FaJoYa06b}, which combines features of the \Struct{} and the
\IJ{} interfaces and targets applications with meshes composed of a set of structured
subgrids, e.g, block-structured, overset, and structured adaptive mesh refinement
grids. The \SStruct{} interface also supports multi-variable PDEs with degrees of freedom
lying in the center, corners, edges or faces of cells composing logically rectangular
boxes. From a computational perspective, these variable types are associated with boxes
that are shifted by different offset values. Thus, we consider only cell-centered problems
here for ease of exposition. The current CPU implementation of SSAMG cannot deal with
problems involving multiple variable types yet; however, the mathematical algorithm of
SSAMG expands to such general cases.

\begin{figure}[!ht]
\centering \includegraphics[scale=0.6]{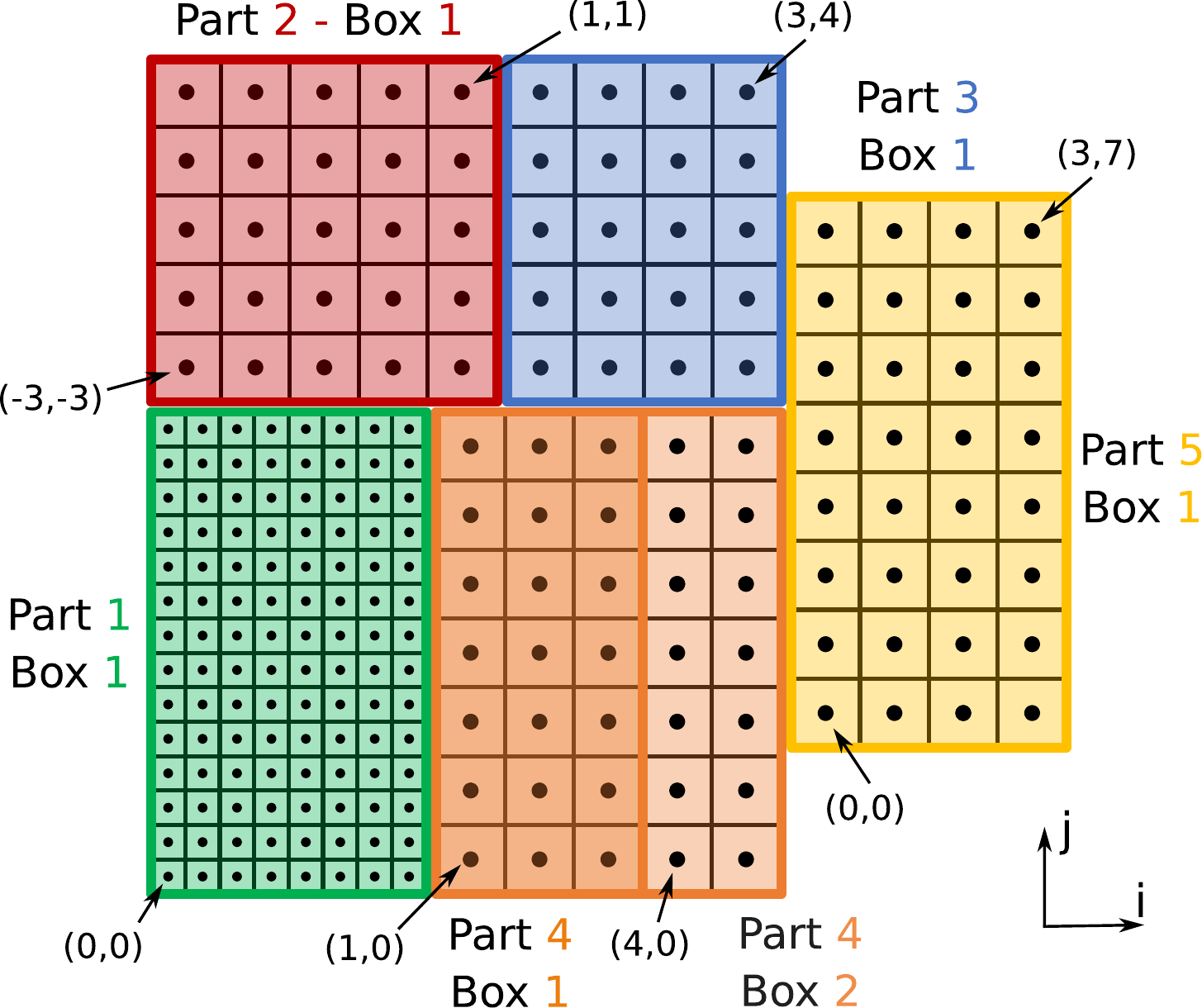}
\caption{A semi-structured grid composed of five parts. Part 4 (orange) consists of two
  boxes, while the others consist of just a single box. Furthermore, Part 1 (green) has a
  refinement factor of two with respect to the other parts. The pairs $(x, y)$ denote cell
  coordinates in the $i$ and $j$ topological directions, respectively. Note that the
  indices of lower-left cells for each part are independent, since the grid parts live in
  different index spaces.}
\label{fig:ssgrid}
\end{figure}

There are five fundamental components required to define a linear system in the \SStruct{}
interface: a grid, stencils, a graph, a matrix, and a vector. The grid is composed of
$n_p$ structured parts with independent index spaces and grid spacing. Each part is formed
topologically by a group of boxes, which are a collection of cell-centered indices,
described by their ``lower'' and ``upper'' corners. Figure \ref{fig:ssgrid} shows an
example of a problem geometry that can be represented by this interface.
Stencils are used to define connections between neighboring grid cells of the same part,
e.g., a typical five-point stencil would connect a generic grid cell to itself and its
immediate neighbors to the west, east, south, and north. The graph describes how
individual parts are connected, see Figure \ref{fig:part-boundaries} for an example. We
have now the components to define a semi-structured matrix $\mat{A} = \mat{S} + \mat{U}$,
which consists of structured and unstructured components, respectively. $\mat{S}$ contains
coefficients that are associated with stencil entries. These can be variable coefficients
for each stencil entry in each cell within a part or can be set to just a single value if
the stencil entry is constant across the part. $\mat{U}$ is stored in \ParCSR{} format and
contains the connections between parts. Lastly, a semi-structured vector describes an
array of values associated with the cells of a semi-structured grid.

\section{Semi-structured algebraic multigrid (SSAMG)}
\label{sec:ssamg}

In the \hypre{} package, there is currently a single native preconditioner for solving
problems with multiple parts through the \SStruct{} interface, which is a block Jacobi
method named \Split{}.  It uses one V-cycle of a structured multigrid solver as an
approximation to the inverse of the structured part of $A$.  This method has limited
robustness since it considers only structured intra-grid couplings in a part to build an
approximation of $\matI{A}$. In this paper, we present a new solver option for the
\SStruct{} interface that computes a multigrid hierarchy taking into account inter-part
couplings. This method is called SSAMG (Semi-Structured Algebraic MultiGrid). It is
currently available in the \texttt{recmat} branch of \hypre{}. This section defines
coarsening, interpolation, and relaxation for SSAMG (subsections~\ref{subsec:coarsening},
\ref{subsec:interp}, and \ref{subsec:relax}, respectively). It also describes how coarse
level operators are constructed (subsection~\ref{subsec:coarseOp}) and discusses a
strategy for improving the method's efficiency at coarse levels
(subsection~\ref{subsec:hybrid}).

\subsection{Coarsening}
\label{subsec:coarsening}

As in PFMG \cite{AsFa96}, we employ semi-coarsening in SSAMG. The coarsening directions
are determined independently for each part of the \SStructGrid{} to allow better treatment
of problems with different anisotropies among the parts. The idea of semi-coarsening is to
coarsen in a single direction of strong coupling such that every other perpendicular
line/plane (2D/3D) forms the new coarse level. For an illustration, see Figure
\ref{fig:part-boundaries}, where coarse points are depicted as solid circles.

In the original PFMG algorithm, the coarsening direction was chosen to be the dimension
with smallest grid spacing. This option is still available in \hypre{} by allowing users
to provide an initial $n_d$-dimensional array of ``representative grid spacings'' that are
only used for coarsening. However, both PFMG and SSAMG can also compute such an array
directly from the matrix coefficients.  In SSAMG, this is done separately for each part,
leading to a matrix $\mat{W} \in \mathbb{R}^{n_p \times n_d}$, where $n_p$ and $n_d$
denote the number of parts and problem dimensions. Here, element $\mat{W}_{pd}$ is
heuristically thought of as a grid spacing for dimension $d$ of part $p$, and hence a
small value indicates strong coupling.

To describe the computation of $\mat{W}$ in part $p$, consider the two-dimensional
nine-point stencil in Figure~\ref{fig:interp}c and assume that $\mat{A}_{C} > 0$ (simple
sign adjustments can be made if $\mat{A}_{C} < 0$). The algorithm extends naturally to
three dimensions. Note also that both PFMG and SSAMG are currently restricted to stencils
that are contained within this nine-point stencil (27-point in 3D). The algorithm proceeds
by first reducing the nine-point matrix to a single five-point stencil through an
averaging process, then computing the (negative) sum of the resulting off-diagonal
coefficients in each dimension. That is, for the $i$-direction ($d = 1$), we compute
\begin{equation}
  c_1 = \sum_{(i,j)} - (\mat{A}_{SW} + \mat{A}_{W} + \mat{A}_{NW}) -
                     (\mat{A}_{SE} + \mat{A}_{E} + \mat{A}_{NE}),
\end{equation}
where the stencil coefficients are understood to vary at each point $(i,j)$ in the grid.
Here the left and right parenthetical sums contribute to the ``west'' and ``east''
coefficients of the five-point stencil. The computation is analogous for the
$j$-direction.  From this, we define
\begin{equation}
  \mat{W}_{pd} = \sqrt{ \dfrac{\max\limits_{\substack{0 \leq i < n_d}} c_i}{c_d}},
\end{equation}
based on the heuristic that the five-point stencil coefficients are inversely proportional
to the square of the grid spacing.

With $\mat{W}$ in hand, the semi-coarsening directions for each level and part are
computed as described in Algorithm \ref{algo:coarsening}.  The algorithm starts by
computing a bounding box\footnote{Given a set of boxes, a bounding box is defined by the
  cells with minimum index (lower corner) and maximum index (upper corner) over the entire
  set.} around the grid in each part, then loops through the grid levels from finest
(level~$0$) to coarsest (level~$n_l$). For a given grid level~$l$ and part~$p$, the
coarsening direction~$d^{\star}$ is set to be the one with minimum\footnote{In the case of
  two or more directions sharing the same value of $\mat{W}_{pd}$, as in an isotropic
  scenario, we set $d^{\star}$ to the one with smallest index.} value in $\mat{W}_{p}$
(line~8). Then, the bounding box for part~$p$ is coarsened by a factor of two in
direction~$d^{\star}$ (line~9) and $\mat{W}_{p,d^{\star}}$ is updated to reflect the
coarser ``grid spacing'' on the next grid level (line~10). If the bounding box is too
small, no coarsening is done (line~7) and that part becomes inactive. The coarsest grid
level $n_l$ is the first level with total semi-structured grid size less than a given
maximum size $s_{max}$, unless this exceeds the specified maximum number of levels
$l_{max}$.

\begin{algorithm}[ht!]
\caption{\bf SSAMG coarsening}
\label{algo:coarsening}
\begin{algorithmic}[1]
\Procedure{SSAMGCoarsen}{$\mat{W}$} \For{$p = 1,n_{p}$} \State Compute part bounding boxes
$bbox_{p}$ \EndFor \For{$l = 1,n_{l}$} \For{$p = 1,n_{p}$} \If {$\vol{bbox_{p}} > 1$}
\State $d^{\star} = \argminD{\mat{W}_{pd}}$ \State Coarsen $bbox_{p}$ in direction
$d^{\star}$ by a factor of 2 \State $\mat{W}_{pd^{\star}} = 2 * \mat{W}_{pd^{\star}}$
\EndIf \EndFor \EndFor \EndProcedure
\end{algorithmic}
\end{algorithm}

\subsection{Interpolation}
\label{subsec:interp}

A key ingredient in multigrid methods is the interpolation (or prolongation) operator
$\mat{P}$, the matrix that transfers information from a coarse level in the grid hierarchy
to the next finer grid.  The restriction operator $\mat{R}$ moves information from a given
level to the next coarser grid. For a numerically scalable method, error modes that are
not efficiently reduced by relaxation should be captured in the range of $\mat{P}$, so
they can be reduced on coarser levels \cite{BrHeMc00}.

In SSAMG, we employ a structured operator-based method for constructing prolongation
similar to the method used in \cite{AsFa96}.  It is ``structured'' because $\mat{P}$ is
composed of only a structured component; interpolation is only done within a part, not
between them.  It is ``operator-based'' because the coefficients are algebraically
computed from $\mat{S}$ and are able to capture heterogeneity and anisotropy.  In
\hypre{}, $\mat{P}$ is a rectangular matrix defined by two grids (domain and range), a
stencil, and corresponding stencil coefficients. In the case of $\mat{P}$, the domain grid
is the coarse grid and the range grid is the fine grid.  Since SSAMG uses semi-coarsening,
the stencil for interpolation consists of three coefficients that are computed by
collapsing the stencil of $\mat{A}$, a common procedure for defining interpolation in
algebraic multigrid methods.

\begin{figure}[!hbtp]
\centering \includegraphics[scale=0.7]{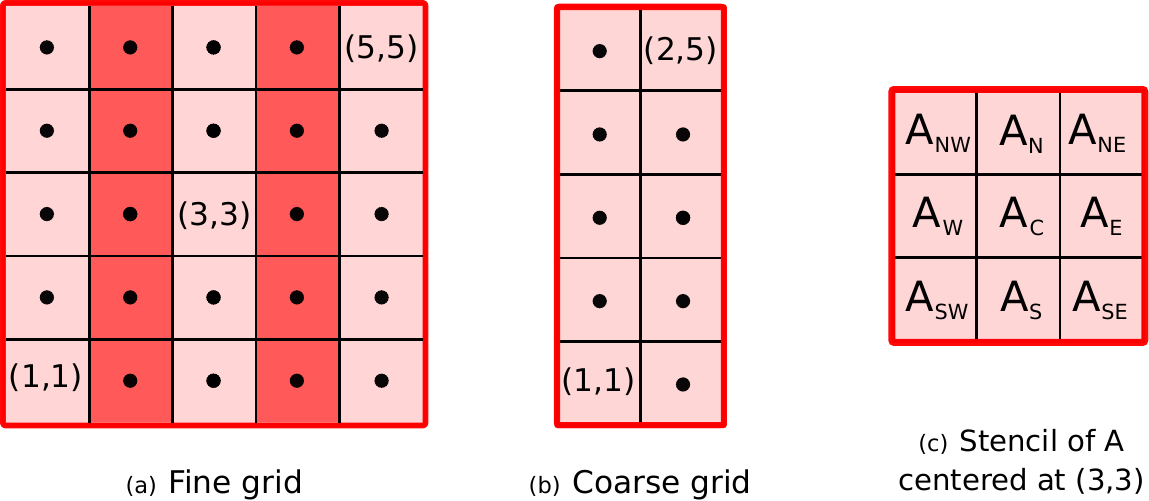}
\caption{(a) and (b) show one example of fine and coarse grids, respectively, also known
  as range and domain grids for the purpose of prolongation. Coarsening is done in the
  $i$-direction, as depicted by the darker cells in the fine grid. (c) shows the stencil
  coefficients of $\mat{A}$ relative to the grid point $(3,3)$ from the fine grid. Stencil
  coefficients for a given grid point can be viewed as the nonzero coefficients of its
  respective row in a sparse matrix.}
\label{fig:interp}
\end{figure}

To exemplify how $\mat{P}$ is computed, consider the solution of the Poisson equation on a
cell-centered grid (Figure \ref{fig:interp}a) formed by a single part and box. Dirichlet
boundary conditions are used and discretization is performed via the finite difference
method with a nine-point stencil (Figure \ref{fig:interp}c).
Assume that coarsening is in the $i$-direction by selecting fine grid cells with even
$i$-coordinate index (depicted in darker red) and renumbering them on the coarse grid as
shown in Figure \ref{fig:interp}b. The prolongation operator connects fine grid cells to
their neighboring coarse grid cells with the following stencil (see~\cite{EnFaYa17} for
more discussion of stencil notation)
\begin{equation*}
  \mat{P} \sim \big[ \mat{P}_W \quad 1 \quad \mat{P}_E \big]_c =
  \begin{bmatrix}
    \mat{P}_W & * & \mat{P}_E
  \end{bmatrix}^{r_1}_c
  ~\oplus~
  \begin{bmatrix}
    * & 1 & *
  \end{bmatrix}^{r_2}_c,
\end{equation*}
where
\begin{equation}
\mat{P}_W = \dfrac{\mat{A}_{SW} + \mat{A}_{W} + \mat{A}_{NW}} {\mat{A}_{S} + \mat{A}_{C} +
  \mat{A}_{N}}, \; \text{and} \quad \mat{P}_E = \dfrac{\mat{A}_{SE} + \mat{A}_{E} +
  \mat{A}_{NE}} {\mat{A}_{S} + \mat{A}_{C} + \mat{A}_{N}}.
\label{eq:prolEx}
\end{equation}
Here, $r_1$ denotes the range subgrid given by the light-colored cells in Figure
\ref{fig:interp}a, and $r_2$ denotes the subgrid given by the dark-colored cells.  For a
fine-grid cell such as $(3,3)$ in Figure~\ref{fig:interp}, interpolation applies the
weights $\mat{P}_W$ and $\mat{P}_E$ to the coarse-grid unknowns associated with cells
$(2,3)$ and $(4,3)$ in the fine-grid indexing, or $(1,3)$ and $(2,3)$ in the coarse-grid
indexing.  For a fine-grid cell such as $(2,3)$, interpolation applies weight~1 to the
corresponding coarse-grid unknown.

When one of the stencil entries crosses a part boundary that is not a physical boundary,
we set the coefficient associated with it to zero and update the coefficient for the
opposite stencil entry so that the vector of ones is contained in the range of the
prolongation operator. Although this gives a lower order interpolation along part
boundaries, it limits stencil growth and makes the computation of coarse level matrices
cheaper, see section \ref{subsec:coarseOp}. It also assures that the near kernel of
$\mat{A}$ is well interpolated between subsequent levels.

Another component needed in a multigrid method is the restriction operator, which maps
information from fine to coarse levels. SSAMG follows the Galerkin approach, where
restriction is defined as the transpose of prolongation $(\mat{R} = \matT{P})$.

\subsection{Coarse level operator}
\label{subsec:coarseOp}

The coarse level operator $\mat{A}_c$ in SSAMG is computed via the Galerkin product
$\matT{P} \mat{A} \mat{P}$. Since the prolongation matrix consists only of the structured
component, the triple-matrix product can be rewritten as
\begin{equation}
  \mat{A}_c = \matT{P}\mat{S}\mat{P} + \matT{P}\mat{U}\mat{P},
\label{eq:galerkin}
\end{equation}
where the first term on the right-hand side is the structured component of $\mat{A}_c$,
and the second its unstructured component. Note that the last term involves the
multiplication of matrices of different types, which we resolve by converting one matrix
type to the other. Since it is generally not possible to represent a \ParCSR{} matrix in
structured format, we convert the structured matrix $\mat{P}$ to the \ParCSR{}
format. However, we consider only the entries of $\mat{P}$ that are actually involved in
the triple-matrix multiplication $\matT{P}\mat{U}\mat{P}$ to decrease the computational
cost of the conversion process.

If we examine the new stencil size for $A_c$, we note that the use of the two-point
interpolation operator limits stencil growth. For example, in the case of a 2D five-point
stencil at the finest level, the maximum stencil size on coarse levels is nine, and for a
3D seven-point stencil at the finest level, the maximum stencil size on coarse levels
is~27.

\begin{figure}
\centering
\includegraphics[scale=0.6]{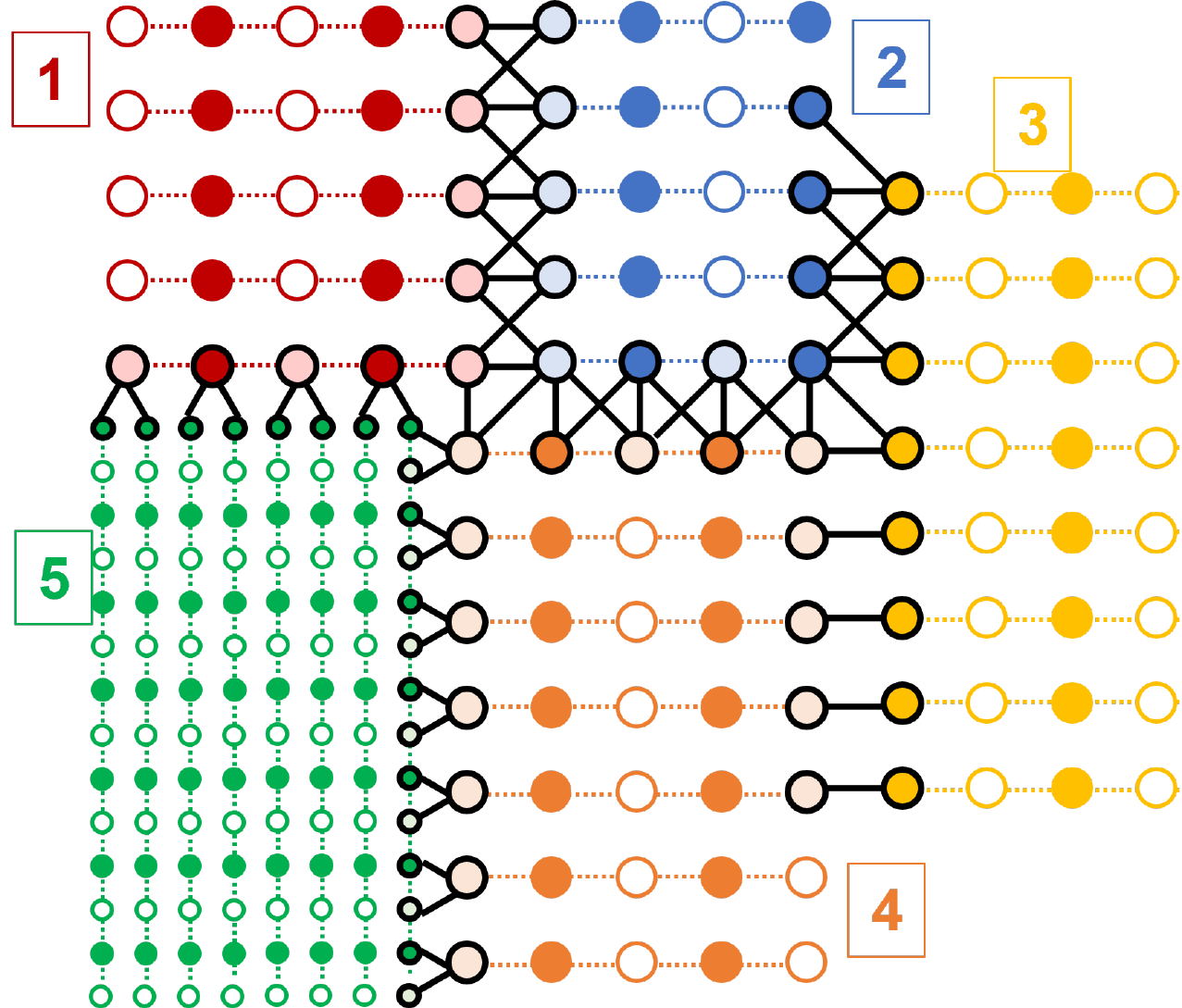}
\caption{Example of a graph of the matrix $\mat{U}$ and graph of matrix $\mat{P}$ derived
  from the semi-structured grid shown in Figure \ref{fig:ssgrid}. The graph of $\mat{U}$
  is depicted by black-solid edges. The graph of $\mat{P}$ consists of five unconnected
  subgraphs illustrated by the dotted multicolored lines. Lastly, the boundary points are
  depicted by black-rimmed circles.}
\label{fig:part-boundaries}
\end{figure}

We prove here that under certain conditions, the unstructured portion of the coarse grid
operator stays restricted to the part boundaries and does not grow into the interior of
the parts. Note that we define a part boundary $\delta \Omega_i$ here as the set of points
in a part $\Omega_i$ that are connected to neighboring parts in the graph of the matrix
$\mat{U}$. For an illustration, see the black-rimmed points in Figure
\ref{fig:part-boundaries}. Figure \ref{fig:part-boundaries} shows the graph of $\mat{P}$
for the semi-structured grid in Figure \ref{fig:ssgrid} and an example of a graph for the
unstructured matrix $\mat{U}$.

\begin{theorem}
We make the following assumptions:
\begin{itemize}
  \item The grid $\Omega$ consists of $k$ parts: $\Omega = \Omega_1 \cup ... \cup
    \Omega_k$, where $\Omega_i \cap \Omega_j = \emptyset$.
  \item The grid has been coarsened using semi-coarsening.
  \item The interpolation $\mat{P}$ interpolates fine points using at most two adjacent
    coarse points aligned with the fine points and maps coarse points onto themselves.
  \item The graph of the unstructured matrix $\mat{U}$ contains only connections between
    boundary points, i.e., $u_{i,j}=0$ if $i \in \Omega_m \setminus \delta \Omega_m,
    m=1,...,k$, or $j \in \Omega_n \setminus \delta\Omega_n, n=1,...,k$, and there are no
    connections within a part, i.e., $u_{i,j}=0$ for $i,j \in \Omega_m, m=1,...,k.$
\end{itemize}
Then the graph of the unstructured part $\mat{U}_{c} = \matT{P}\mat{U}\mat{P}$ also
contains only connections between boundary points, i.e., $u^c_{i,j}=0$ if $i \in
\Omega^c_m \setminus \delta \Omega^c_m, m=1,...,k$, or $j \in \Omega^c_n \setminus
\delta\Omega^c_n, n=1,...,k$, and there are no connections within a part, i.e.,
$u^c_{i,j}=0$ for $i,j \in \Omega^c_m, m=1,...,k.$.
\end{theorem}

\begin{proof}
Since we want to examine how boundary parts are handled, we reorder the interpolation
matrix $\mat{P}$ and the unstructured part $\mat{U}$, so that all interior points are
first followed by all boundary points. The matrices $\mat{P}$ and $\mat{U}$ are then
defined as follows:
\begin{equation}
\mat{P} =
\left(
\begin{array}{cc}
  \mat{P}^I   & \mat{P}^{IB} \\
  \mat{P}^{BI} & \mat{P}^B   \\
\end{array}
\right),
~~~~
\mat{U} =
\left(
\begin{array}{cc}
0 & 0 \\
0 & \mat{U}^B \\
\end{array}
\right).
\end{equation}
Note that while $U^B$ maps $\delta\Omega$ onto $\delta\Omega$, $P^B$ maps $\delta\Omega_c$
onto $\delta\Omega$. Thus, in the extreme case that all boundary points are fine points,
$P^{IB}$ and $P^B$ do not exist. Then, the coarse unstructured part is given as follows:
\begin{equation}
\mat{U}_{c} = \matT{P}\mat{U}\mat{P} =
\left(
\begin{array}{cc}
  (P^{BI})^T \mat{U}^B P^{BI} & (P^{BI})^T \mat{U}^B P^B \\
  (P^B)^T \mat{U}^B P^{BI} & (P^B)^T \mat{U}^B P^{B} \\
\end{array}
\right).
\end{equation}
It is clear already that there is no longer a connection to $\mat{P}^I$ and
$\mat{P^{IB}}$, eliminating many potential connections to interior points; however, we
still need to investigate further the influence of $\mat{P}^{BI}$ and $\mat{P}^{B}$.

Since $P^{BI}$, $P^B$, and $\mat{U}^B$ are still very complex due to their dependence on
$k$ parts, we further define them as follows using the fact that $\mat{P}$ is defined only
on the structured parts and $\mat{U}$ only connects boundary points of neighboring parts.
\begin{equation}
\mat{P}^x =
\left(
\begin{array}{cccc}
P^x_1 &  &  & \\
 & P^x_2 &  & \\
 & & \ddots & \\
 &  & & P^x_k \\
\end{array}\right),
~~~~
\mat{U}^B =
\left(
\begin{array}{cccc}
0 & U^B_{1,2} & ... & U^B_{1,k} \\
U^B_{2,1} & 0 & \ddots & \vdots  \\
\vdots & \ddots & \ddots & U^B_{k-1,k} \\
U^B_{k,1} & ... & U^B_{k,k-1}  & 0 \\
\end{array}\right).
\end{equation}
Note that while $U_{ij}^B$ maps $\delta \Omega_i$ to $\delta \Omega_j$, only the
coefficients corresponding to edges in the graph of $U$ that connect points in $\delta
\Omega_i$ to $\delta \Omega_j$ are nonzero, all other coefficients are zero.  Then,
$(P^x)^T\mat{U}^BP^{y}$, where ``$x$'' and ``$y$" can stand for ``$BI$" as well as ``$B$",
is given by
\begin{equation}
\left(
\begin{array}{cccc}
	0 & (P^x_1)^TU^B_{1,2} P^{y}_2& ... & (P^x_1)^TU^B_{1,k} P^{y}_k\\
	(P^x_2)^TU^B_{2,1}P^{y}_1 & 0 & \ddots & \vdots  \\
	\vdots & \ddots & \ddots & (P^x_{k-1})^TU^B_{k-1,k} P^{y}_k\\
	(P^x_k)^TU^B_{k,1} P^{y}_1& ... & (P^x_k)^TU^B_{k,k-1} P^{y}_{k-1} & 0 \\
\end{array}
\right).
\label{eqn:PiUPj}
\end{equation}
This allows us to just focus on the submatrices $(P^x_i)^TU^B_{ij}P^{y}_j$.  There are
three potential scenarios that can occur at the boundary between two parts due to our use
of semi-coarsening and a simple two-point interpolation (Figure
\ref{fig:part-boundaries}):
\begin{itemize}
  \item all boundary points are coarse points as shown at the right boundary of part 2 and
    the left boundary of part 3;
  \item all boundary points are fine points as at the right boundary of part 1 and 4;
  \item the boundary points are alternating coarse and fine points as illustrated at the
    right boundary of part 5.
\end{itemize}

Let us define $P_i^x|_{\delta \Omega_{ij}}$ as the matrix that consists of the rows of
$P_i^x$ that correspond to all boundary points in $\delta \Omega_i$ that are connected to
boundary points in $\delta \Omega_j$.  If all points are coarse points, $P_i^B|_{\delta
  \Omega_{ij}}=I$ and $P_i^{BI}|_{\delta \Omega_{ij}} = 0$, since there are no connections
from the boundary to the interior for $P_i^{BI}|_{\delta \Omega_{ij}}$.  If all points are
fine points, $P_i^B|_{\delta \Omega_{ij}}$ does not exist, and $P_i^{BI}|_{\delta
  \Omega_{ij}}$ is a matrix with at most one nonzero element per row.  Since coarse points
in $\Omega_i$ adjacent to the fine boundary points in $\delta \Omega_i$ become boundary
points of $\Omega_i^c$, e.g., see right boundary of part 1 or left and right boundaries of
part 4, all nonzero elements in $P_i^{BI}|_{\delta \Omega_{ij}}$ are associated with a
column belonging to $\delta \Omega_i$.  In the case of alternating fine and coarse points,
$P_i^{BI}|_{\delta \Omega_{ij}}=0$, since there are no connections from the boundary to
the interior, and $P_i^B|_{\delta \Omega_{ij}}$ is a matrix with at most two nonzeros in
the $j$-th and $k$-th columns, where $j$ and $k$ are elements of $\delta \Omega_i^c$.
Recall that all columns in $U_{ij}$ belonging to points outside of $\delta \Omega_j$ and
all rows belonging to points outside of $\delta \Omega_i$ are zero. Based on this and the
previous observations it is clear that if all points are coarse or we are dealing with
alternating fine and coarse points, the submatrices in \ref{eqn:PiUPj} that involve
$P_i^{BI}$ will be 0, since $P_i^{BI}|_{\delta \Omega_{ij}}=0$ and $P_j^{BI}|_{\delta
  \Omega_{ji}}=0$. Any additional nonzero coefficients in $P_i^{BI}$ or $P_j^{BI}$ due to
boundary points next to other parts will be canceled out in the matrix product.  It is
also clear, since the columns of $P_i^B$ pertain only to points in $\delta \Omega_i^c$,
that the graph of the product $(P_i^B)^TU_{ij}P_j^B$ only contains connections of points
of $\delta \Omega_i^c$ to points of $\delta \Omega_j^c$ and none to the interior or to
itself.

Let us further investigate the case where all boundary points are fine points.  We first
consider $(P_i^{BI})^TU_{ij}P_j^{BI}$. Since we have already shown that $P_i^{BI}|_{\delta
  \Omega_{ij}}=0$ for boundaries with coarse or alternating points leading to zero triple
products in \ref{eqn:PiUPj}, we can ignore these scenarios and assume that for both
$P_i^{BI}$ and $P_j^{BI}$ the boundary points adjacent to each other are fine points.
Each row of $P_i^{BI}|_{\delta \Omega_{ij}}$ has at most one nonzero element in the column
corresponding to the interior coarse point connected to the fine boundary point. This
interior point is also an element in $\delta \Omega_i^c$.  Therefore the graph of the
product $(P_i^{BI})^TU_{ij}P_j^{BI}$ only contains connections of points of $\delta
\Omega_i^c$ to points of $\delta \Omega_j^c$ and none to the interior or to itself.
Finally, this statement also holds for the triple products $(P_i^B)^TU_{ij}P_j^{BI}$ and
$(P_i^{BI})^TU_{ij}P_j^B$ using the same arguments as above.

Note that the number of nonzero coefficients in $\mat{U}_c$ can still be larger than those
in $\mat{U}$, however the growth only occurs along part boundaries.
\end{proof}

\subsection{Relaxation}
\label{subsec:relax}

Relaxation, or smoothing, is an important element of multigrid whose task is to eliminate
high frequency error components from the solution vector $\vec{x}$. The relaxation process
at step $k > 0$ can be described via the generic formula:
\begin{equation}
  \vec{x}_{k} = \vec{x}_{k-1} + \omega \matI{M} \left( \vec{b} - \mat{A} \vec{x}_{k-1}
  \right),
\label{eq:relax0}
\end{equation}
where $\matI{M}$ is the smoother operator and $\omega$ is the relaxation weight. In SSAMG,
we provide two pointwise relaxation schemes. The first one is weighted Jacobi, in which
$\matI{M} = \matI{D}$, with $\mat{D}$ being the diagonal of $\mat{A}$. Moreover, $\omega$
varies for each multigrid level and semi-structured part as a function of the grid-spacing
metric $\mat{W}$:

\begin{equation}
  \omega_{p} = \dfrac{2}{3 - \beta_{p}/\alpha_{p}},
\label{eq:relax1}
\end{equation}
where
\begin{equation}
  \alpha_{p} = \sum^{n_d}_{d = 0} \dfrac{1}{\mat{W}^2_{pd}} \quad \text{and} \quad
  \beta_{p} = \sum^{n_d}_{\mathclap{\substack{d \,=\, 0, \\d \neq d^{\star}}}}
  \dfrac{1}{\mat{W}^2_{pd}}.
\label{eq:alpha-beta}
\end{equation}
The ratio $\beta_p/\alpha_p$ adjusts the relaxation weight to more closely approximate the
optimal weight for isotropic problems in different dimensions. To see how this works,
consider as an example a highly-anisotropic 3D problem that is nearly decoupled in the
$k$-direction and isotropic in~$i$ and~$j$. Because of the severe anisotropy, the problem
is effectively~2D, so the optimal relaxation weight is~$4/5$. Since our coarsening
algorithm will only coarsen in either directions~$i$ or~$j$, we get
$\beta_p/\alpha_p=1/2$, and $\omega_p=4/5$ as desired.

The second relaxation method supported by SSAMG is L1-Jacobi. This method is similar to
the previous one, in the sense that a diagonal matrix is used to construct the smoother
operator; however, here, the $i$-th diagonal element of $\mat{M}$ equals the L1-norm of
the $i$-th row of $A$:
\begin{equation*}
  \mat{M}_{ii} = \sum\limits^{N}_{j = 0} \left|\mat{A}_{ij} \right|.
\label{eq:relax2}
\end{equation*}
This form leads to guaranteed convergence when $A$ is positive definite, i.e., the error
propagation operator $\mat{E} = \mat{I} - \matI{M} \mat{A}$ has a spectral radius smaller
than one. We refer to \cite{BaFaKo11} for more details. This option tends to give slower
convergence than weighted Jacobi; however, a user-defined relaxation factor in the range
$(1, {2}/{\lambda_{max}(\matI{M} \mat{A})})$ ($\lambda_{max}$ is the maximum eigenvalue)
can be used to improve convergence.

To reduce the computational cost of a multigrid cycle within SSAMG, we also provide an
option to turn off relaxation on certain multigrid levels in isotropic (or partially
isotropic) scenarios. We call this option ``skip'', and it has the action of mimicking
full-coarsening. For a given part, SSAMG's algorithm checks if the coarsening directions
for levels ``$l$'' and ``$l - n_d$'' match $\left(d^{\star}_{l} = d^{\star}_{l -
  n_d}\right)$. If yes, then relaxation is turned off (skipped) at level $l$.

\subsection{Hybrid approach}
\label{subsec:hybrid}

Since SSAMG uses semi-coarsening, the coarsening ratio between the number of variables on
subsequent grids is~2. In classical algebraic multigrid, this value tends to be larger,
especially when aggressive coarsening strategies are applied. This leads to the creation
of more levels in the multigrid hierarchy of SSAMG when compared to BoomerAMG. Since the
performance benefits of exploiting structure decreases on coarser grid levels, we provide
an option to transition to an unstructured multigrid hierarchy at a certain level or
coarse problem size chosen by the user. This is done by converting the matrix type from
\SStructMatrix{} to \ParCSRMatrix{} at the transition level. The rest of the multigrid
hierarchy is set up using BoomerAMG configured with the default options used in
\hypre{}. With a properly chosen transition level, this hybrid approach can improve
performance. In the non-hybrid case, SSAMG employs one sweep of the same relaxation method
used in previous levels.

\section{Numerical results}
\label{sec:results}

In this section, we investigate convergence and performance of SSAMG when used as a
preconditioner for the conjugate gradient method (PCG). We also compare it to three other
multigrid schemes in \hypre{}, namely PFMG, Split, and BoomerAMG. The first is the
flagship multigrid method for structured problems in \hypre{} based on semi-coarsening
\cite{AsFa96,BrFaJo00}, the second, a block-Jacobi method built on top of the \SStruct{}
interface~\cite{Hypre}, in which blocks are mapped to semi-structured parts, and the last
scheme is \hypre{}'s unstructured algebraic multigrid method~\cite{HeYa2002}. Each of
these preconditioners has multiple setup parameters that affect its performance. For the
comparison made here, we select those leading to the best solution times on CPU
architectures. In addition, we consider four variants of SSAMG in an incremental setting
to demonstrate the effects of different setup options described in the paper. A complete
list of the methods considered here is given below:

\begin{itemize}
  \item \PFMGOpt{}: weighted Jacobi\footnote{This is the default relaxation method of
    PFMG.} smoother and ``skip'' option turned on.
  \item \SStructSplitOpt{}: block-Jacobi method with one V-cycle of PFMG as the inner solver for
    parts.
  \item \BoomerAMGOpt{}: Forward/Backward L1-Gauss-Seidel relaxation \cite{BaFaKo11};
    coarsening via HMIS \cite{StYaHe06} with a strength threshold value of $0.25$;
    modularized option for computing the Galerkin product $\mat{RAP}$; one level (first)
    of aggressive coarsening with multi-pass interpolation \cite{Ya10} and, in the
    following levels, matrix-based extended+i interpolation \cite{LiSjYa21} truncated to a
    maximum of four nonzero coefficients per row.
  \item \SSAMGBase{}: baseline configuration of SSAMG employing weighted L1-Jacobi
    smoother with relaxation factor equal to $3/2$.
  \item \SSAMGSkip{}: above configuration plus the ``skip'' option.
  \item \SSAMGHybrid{}: above configuration plus the ``hybrid'' option for switching to
    BoomerAMG as the coarse solver at the $10^{\text{th}}$ level, which corresponds to
    three steps of full grid refinement in 3D, i.e., $512$ times reduction on the number
    of degrees of freedom (DOFs).
  \item \SSAMGOpt{}: refers to the best SSAMG configuration and employs the same
    parameters as \SSAMGHybrid{} except for switching to BoomerAMG at the $7^{\text{th}}$
    level. This results in six pure SSAMG coarsening levels and reduction factor of $64$
    on the number of DOFs.
\end{itemize}

Every multigrid preconditioner listed above is applied to the residual vector via a single
V(1,1)-cycle. The coarsest grid size is at most 8 in all cases where BoomerAMG is used, it
equals the number of parts for \SSAMGBase{} and \SSAMGSkip{}, and one for \PFMGOpt{}. The
number of levels in the various multigrid hierarchies changes for different problem sizes
and solvers.

We consider four test cases featuring three-dimensional semi-structured grids, different
part distributions, and anisotropy directions. Each semi-structured part is formed by a
box containing $m \times m \times m$ cells. Similarly, in a distributed parallel setting,
each semi-structured part is owned by $p \times p \times p$ unique MPI tasks, meaning that
the global semi-structured grid is formed by replicating the local $m^3$-sized boxes
belonging to each part by $p$ times in each topological direction. This leads to a total
of $n_p p^3$ MPI tasks for $n_p$ parts. We are particularly interested in evaluating weak
scalability of the proposed method for a few tasks up to a range of thousands of MPI
tasks. Thus, we vary the value of $p$ from one to eight with increments of one.

For the results, we report the number of iterations needed for convergence, setup time of
the preconditioner, and solve time of the iterative solver.  All experiments were
performed on Lassen, a cluster at LLNL equipped with two IBM POWER9 processors (totaling
$44$ physical cores) per node. However, we note that up to $32$ cores per node were used
in the numerical experiments to reduce the effect of limited memory bandwidth. Convergence
of the iterative solver is achieved when the L2-norm of the residual vector is less than
$10^{-6} ||\vec{b}||_2$. The linear systems were formed via discretization of the Poisson
equation through the finite difference method, and zero Dirichlet boundary conditions are
used everywhere except for the $k = 0$ boundary, which assumes a value of one. The initial
solution guess passed to PCG is the vector of zeros. The discretization scheme we used
leads to the following seven-point stencil:
\begin{equation}
\mat{A} \sim
\begin{bmatrix}
-\gamma
\end{bmatrix}
\begin{bmatrix}
        & -\beta & \\
-\alpha & 2(\alpha + \beta + \gamma) & -\alpha \\
        & -\beta &
\end{bmatrix}
\begin{bmatrix}
-\gamma
\end{bmatrix}
\label{eq:stencilA}
\end{equation}
where $\alpha$, $\beta$, and $\gamma$ denote the coefficients in the $i$, $j$,
and $k$ topological directions. For the isotropic problems, $\alpha = \beta = \gamma = 1$,
for the anisotropic cases we define their values in section \ref{subsec:cubesAniso}.

\subsection{Test case 1 - cubes side-by-side}
\label{subsec:cubes}

\begin{figure}[!hbtp]
\centering
\includegraphics[trim={0 1.1cm 11cm 0}, clip, scale=0.5]{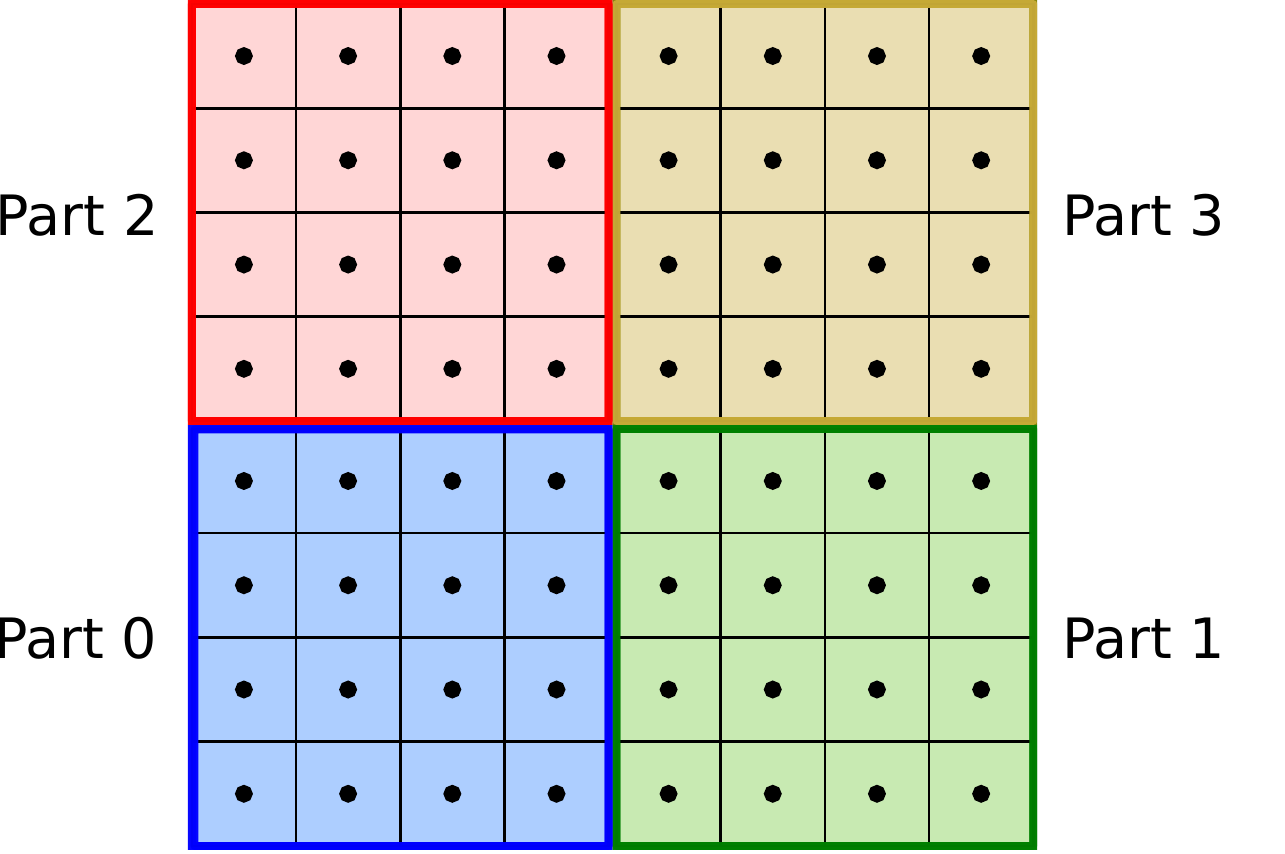}
\caption{$ij$-plane cut of the three-dimensional base grid used for test case
  \ref{subsec:cubes}. There are no adjacent parts in the $k$-direction. Colors denote
  different parts, and the experiments showed in this section are produced by equally
  refining the base grid in all directions.}
\label{fig:cubes4}
\end{figure}

The first test case is made of an isotropic and block-structured three-dimensional domain
composed of four cubes, where each contains the same number of cells and refers to a
different semi-structured part. Figure~\ref{fig:cubes4} shows a plane view of one
particular case with cubes formed by four cells in each direction. Regarding the solver
choices, since PFMG works only for single-part problems, we translated parts into
different boxes in an equivalent structured grid. Note that such a transformation is only
possible due to the simplicity of the current problem geometry and is unattainable in more
general cases such as those described later in sections \ref{subsec:miller} and
\ref{subsec:samr}.

For the numerical experiments, we consider $m = 128$, which gives a local problem size per
part of $2,097,152$ DOFs and a global problem size of $8,388,608$ DOFs for 4 MPI tasks ($p
= 1$). The largest problem we consider here, obtained when $p = 8$, has a global size of
about $4.3$B DOFs.

\begin{figure}[!hbtp]
\centering
\includegraphics[trim={4.8cm 0 4.85cm 0}, clip, width=\textwidth]{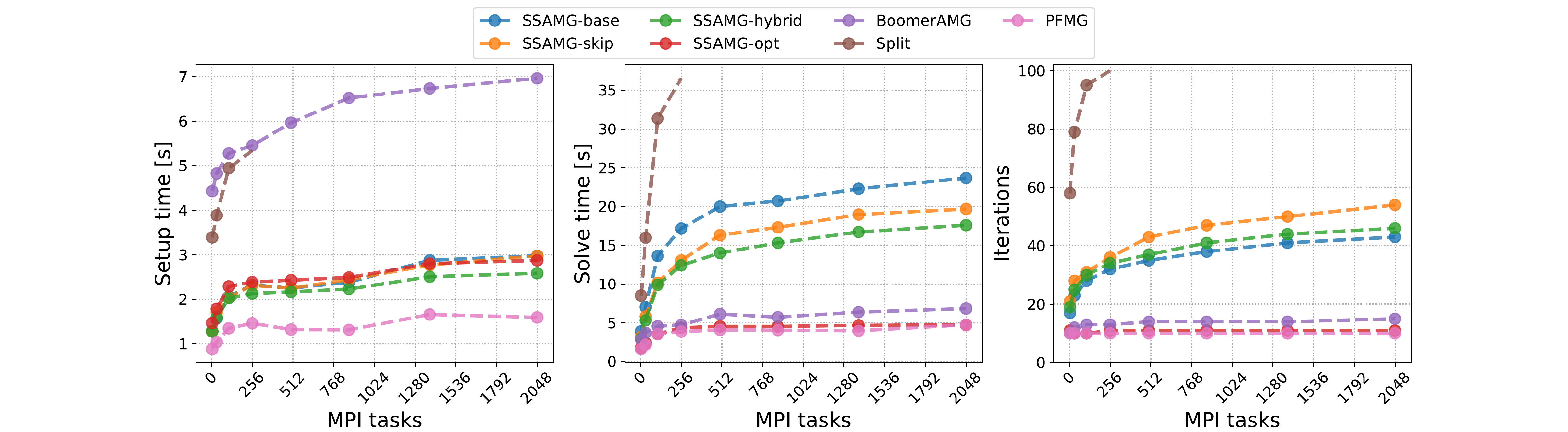}
\caption{Weak scalability results for test case 1. Three metrics are shown in the figure,
  i.e., setup phase times in seconds (left); solve phase times in seconds (middle), and
  number of iterations (right). All curves are plotted with respect to the number of MPI
  tasks, $N_{\text{procs}}$, which varies from $4$ ($p = 1$) up to $2048$ ($p = 8$).}
\label{fig:t0_ref32x32x32}
\end{figure}

Figure~\ref{fig:t0_ref32x32x32} shows weak scalability results for this test
case. Analyzing the iteration counts, \SStructSplitOpt{} is the only method that does not
converge in less than a maximum iteration count of~$100$ for runs larger than $500$M DOFs
($p = 4$). This lack of numerical scalability was already expected since couplings among
parts are not captured in \SStructSplitOpt{}'s multigrid hierarchy.  The best iteration
counts are reached by PFMG, which is natural since this method can take full advantage of
the problem's geometry. Noticeably, the iteration counts of \SSAMGOpt{} follow PFMG
closely, since part boundaries are no longer considered after switching to BoomerAMG on
the coarser levels and the switch is done earlier here than in \SSAMGHybrid{}; the other
SSAMG variants need a higher number of iterations for achieving convergence, since the
interpolation is of lower quality along part boundaries. Lastly, the BoomerAMG
preconditioner shows a modest increase in iteration counts for increasing problem sizes,
and this is common in the context of algebraic multigrid.

Solve times are directly related to iteration counts. Since \SStructSplitOpt{} has a
similar iteration cost to the other methods but takes the largest number of iterations to
converge, it is the slowest option in solution time. For the same reason, the three SSAMG
variants except for \SSAMGOpt{} are slower than the remaining
preconditioners. Still, \SSAMGSkip{} is faster than \SSAMGBase{} despite showing more
iterations because the ``skip''option reduces its iteration cost. The optimal variant
\SSAMGOpt{} is able to beat \BoomerAMGOpt{} by a factor of $1.6$x for $p = 1$ and $2.3$x
for $p = 8$. Moreover, \SSAMGOpt{} shows little performance degradation with respect to
the fastest preconditioner (\PFMGOpt{}).

BoomerAMG is the slowest option when analyzing setup times. This is a result of multiple
reasons, the three most significant ones being:
\begin{itemize}
  \item BoomerAMG employs more elaborate formulas for computing interpolation, which
    require more computation time than the simple two-point scheme used by PFMG and SSAMG;
  \item the triple-matrix product algorithm for computing coarse operators implemented for
    CSR matrices is less efficient than the specialized algorithm employed by \Struct{}
    and \SStruct{} matrices\footnote{We plan to explore this statement with more depth in
      a following communication.};
  \item BoomerAMG's coarsening algorithm involves choosing fine/coarse nodes on the matrix
    graph besides computing a strength of connection matrix. Those steps are not necessary
    for PFMG or SSAMG.
\end{itemize}
This is followed by \SStructSplitOpt{}, which should have setup times close to PFMG, but
due to a limitation of its parallel implementation, the method does not scale well with an
increasing number of parts. On the other hand, all the SSAMG variants show comparable
setup times, up to $2.8$x faster than \BoomerAMGOpt{}. The first two SSAMG variants share
the same setup algorithm, and their lines are superposed. \SSAMGOpt{} has a slightly
slower setup for $p \leq 5$ than \SSAMGBase{}, but for $p > 5$ the setup times of these
two methods match. The fastest SSAMG variant by a factor of $1.2$x with respect to the
others is \SSAMGHybrid{}, and that holds because it generates a multigrid hierarchy with
fewer levels than the non-hybrid SSAMG variants leading to less communication overhead
associated with collective MPI calls. The same argument is true for \SSAMGOpt{}; however,
the benefits of having fewer levels is outweighed by the cost of converting the
\SStructMatrix{} to a \ParCSRMatrix{} in the switching level. Still, \SSAMGOpt{} is $2.3$x
and $3$x faster than \BoomerAMGOpt{} for $p = 1$ and $p = 8$, respectively. Finally,
\PFMGOpt{} yields the best setup times with a speedup of nearly $4.6$x with respect to
BoomerAMG and up to $1.9$x with respect to SSAMG.

We note that PFMG is naturally a better preconditioner for this problem than SSAMG since
it interpolates across part boundaries. However, this test case was significant to show
how close the performance of SSAMG can be to PFMG, and we demonstrated that \SSAMGOpt{} is
not much behind \PFMGOpt{}, besides yielding faster solve and setup times than BoomerAMG.

\subsection{Test case 2 - anisotropic cubes}
\label{subsec:cubesAniso}

\begin{figure}[!hbtp]
\centering
\includegraphics[scale=0.7]{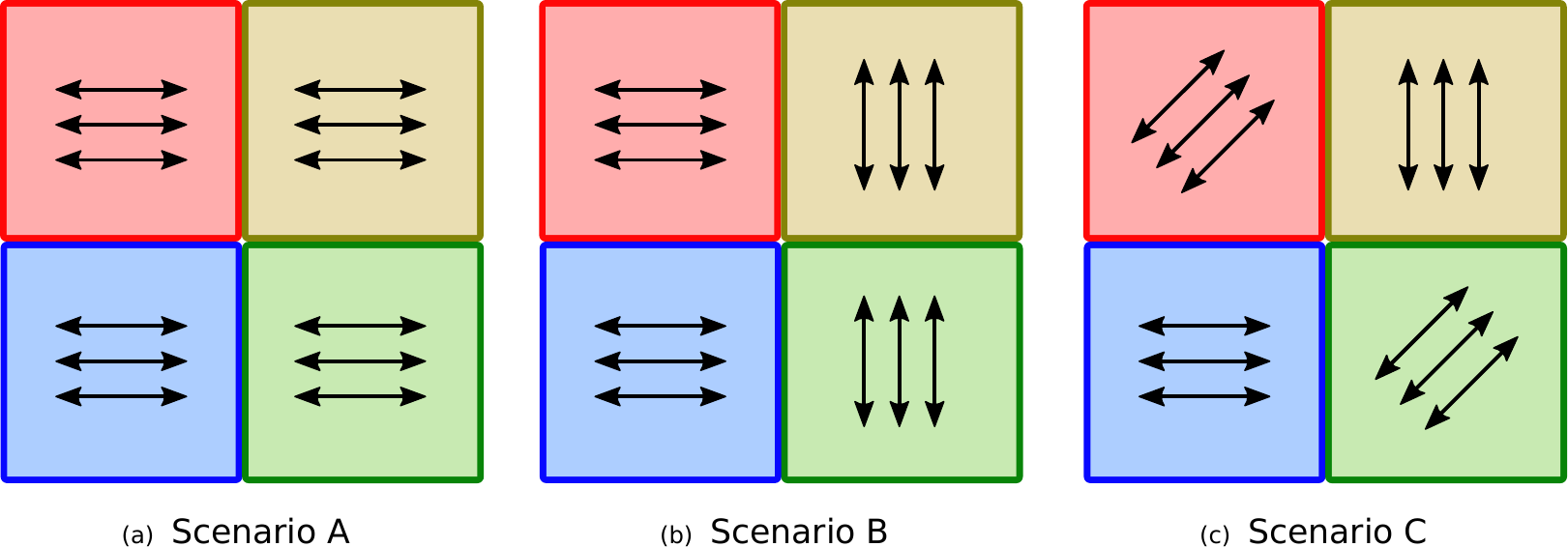}
\caption{$XY$-plane cut of the three-dimensional grids used in test case
  \ref{subsec:cubesAniso}. We consider three anisotropy scenarios. Arrows indicate the
  direction of prevailing anisotropy in each part of the grid, e.g., $i$-direction in
  scenario A. Diagonal arrows in the rightmost case indicate the $k$-direction.}
\label{fig:aniso}
\end{figure}

This test case has the same problem geometry and sizes ($m = 128$) as the previous test
case; however, it employs different stencil coefficients ($\alpha$, $\beta$, and~$\gamma$)
for each part of the grid with the aim of evaluating how anisotropy affects solver
performance. Particularly, we consider three different scenarios (Figure~\ref{fig:aniso})
where the coefficients relative to stencil entries belonging to the direction of strongest
anisotropy for a given part are $100$~times larger than the remaining ones. The directions
of prevailing anisotropy for each scenario are listed below:
\begin{enumerate}[label=(\Alph*)]
  \item ``$i$'' (horizontal) for all semi-structured parts.
  \item ``$i$'' for parts zero and two; ``$j$'' (vertical) for parts one and three.
  \item ``$i$'' for part zero, ``$j$'' for part three, and ``$k$'' (depth) for parts one
    and two.
\end{enumerate}
Regarding the usage of PFMG for this problem, the same transformation mentioned in section
\ref{subsec:cubes} apply here as well.

\begin{figure}[!hbtp]
\centering
\includegraphics[trim={4.8cm 0 4.85cm 0}, clip, width=\textwidth]{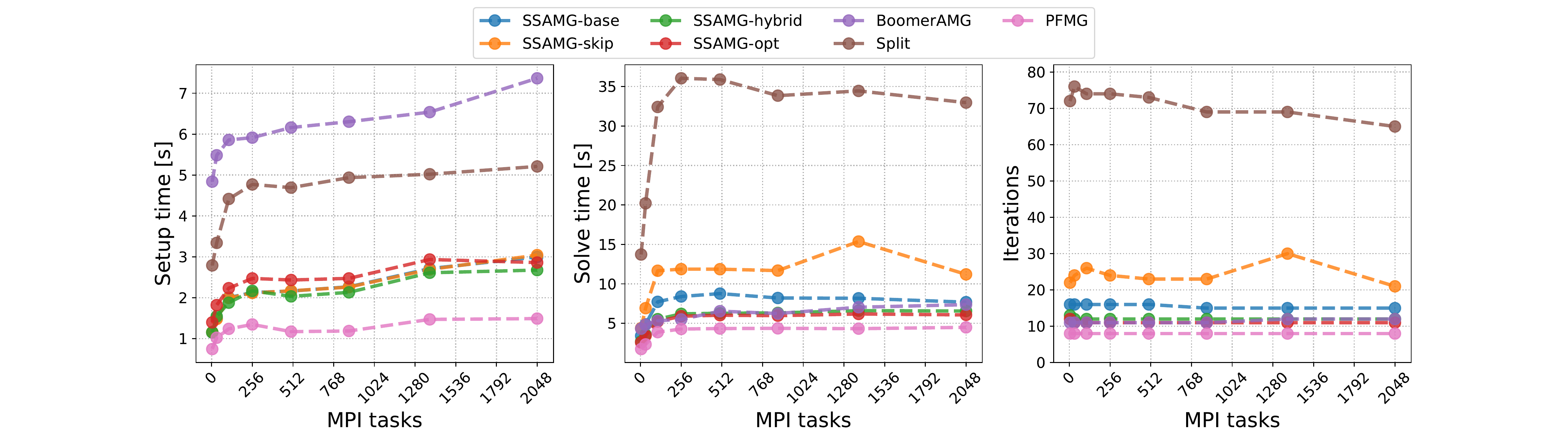}
\caption{Weak scalability results for scenario A of test case 2. Three metrics are shown
  in the figure, i.e., setup phase times in seconds (left); solve phase times in seconds
  (middle), and number of iterations (right). All curves are plotted with respect to the
  number of MPI tasks, $N_{\text{procs}}$, which varies from $4$ ($p = 1$) up to $2048$
  ($p = 8$).}
\label{fig:t1A_ref32x32x32}
\end{figure}

Figure~\ref{fig:t1A_ref32x32x32} shows the results referent to scenario A. The numerical
scalabilities of the different methods look better than in the previous test case. This is
valid especially for the SSAMG variants, because the two-point
interpolation strategy is naturally a good choice for the first few coarsening levels when
anisotropy is present in the matrix coefficients. Such observation also applies to the
less scalable \SStructSplitOpt{} method, explaining the better behavior seen here with
respect to Figure \ref{fig:t0_ref32x32x32}. Again, \PFMGOpt{} uses the least number
of iterations followed closely by \SSAMGOpt{} and \BoomerAMGOpt{}.

Regarding solve times, \SSAMGOpt{} is about $1.3$x faster than \BoomerAMGOpt{} for $p \leq
2$, while for $p > 2$ these methods show similar times. The ``skip'' option of SSAMG is
not beneficial for this case since the solve times of \SSAMGSkip{} are higher than
\SSAMGBase{}. In fact, such an option does not play a significant role in reducing the
solve time compared to isotropic test cases. This is because coarsening happens in the
same direction for the first few levels in anisotropic test cases, and thus relaxation is
skipped only in the later levels of the multigrid hierarchy where the cost per iteration
associated with them is already low compared to the initial levels. Moreover, the omission
of relaxation in coarser levels of the multigrid hierarchy can be detrimental for
convergence in SSAMG, explaining why \SSAMGSkip{} requires more iterations than
\SSAMGBase{}. Following the fact that \PFMGOpt{} is the method that needs fewer iterations
for convergence, it is also the fastest in terms of solution times. For setup times, the
four SSAMG variants show comparable results, and similar conclusions to test case 1 are
valid here. Lastly, the speedups of \SSAMGOpt{} over \BoomerAMGOpt{} are $3.3$x and $2.5$x
for $p = 1$ and $p = 8$, respectively.

\begin{figure}[!hbtp]
\centering
\includegraphics[trim={4.8cm 0 4.85cm 0}, clip, width=\textwidth]{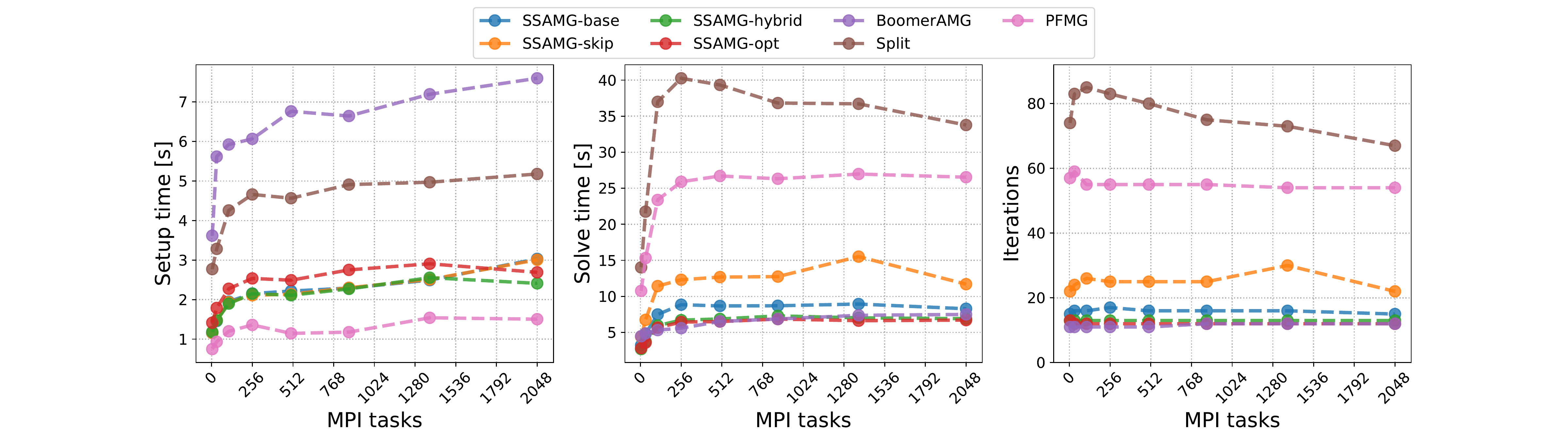}
\caption{Weak scalability results for scenario B of test case 2. Three metrics are shown
  in the figure, i.e., setup phase times in seconds (left); solve phase times in seconds
  (middle), and number of iterations (right). All curves are plotted with respect to the
  number of MPI tasks, $N_{\text{procs}}$, which varies from $4$ ($p = 1$) up to $2048$
  ($p = 8$).}
\label{fig:t1B_ref32x32x32}
\end{figure}

Results for scenario B are shown in Figure~\ref{fig:t1B_ref32x32x32}. The most significant
difference here compared to scenario A are the results for \PFMGOpt{}. Particularly, the
number of iterations is much higher than in the previous cases. This is caused by the fact
that PFMG employs the same coarsening direction everywhere on the grid, and thus it cannot
recognize the different regions of anisotropy as done by SSAMG. This is clearly
sub-optimal since a good coarsening scheme should adapt to the direction of largest
coupling of the matrix coefficients. The larger number of iterations is also reflected in
the solve times of PFMG, which become less favorable than those by SSAMG and
BoomerAMG. Setup times of PFMG continue to be the fastest ones; however, this advantage is
not sufficient to maintain its position of fastest method overall. The comments regarding
the speedups of SSAMG compared to BoomerAMG made for scenario A also apply here.

\begin{figure}[!hbtp]
\centering
\includegraphics[trim={4.8cm 0 4.85cm 0}, clip, width=\textwidth]{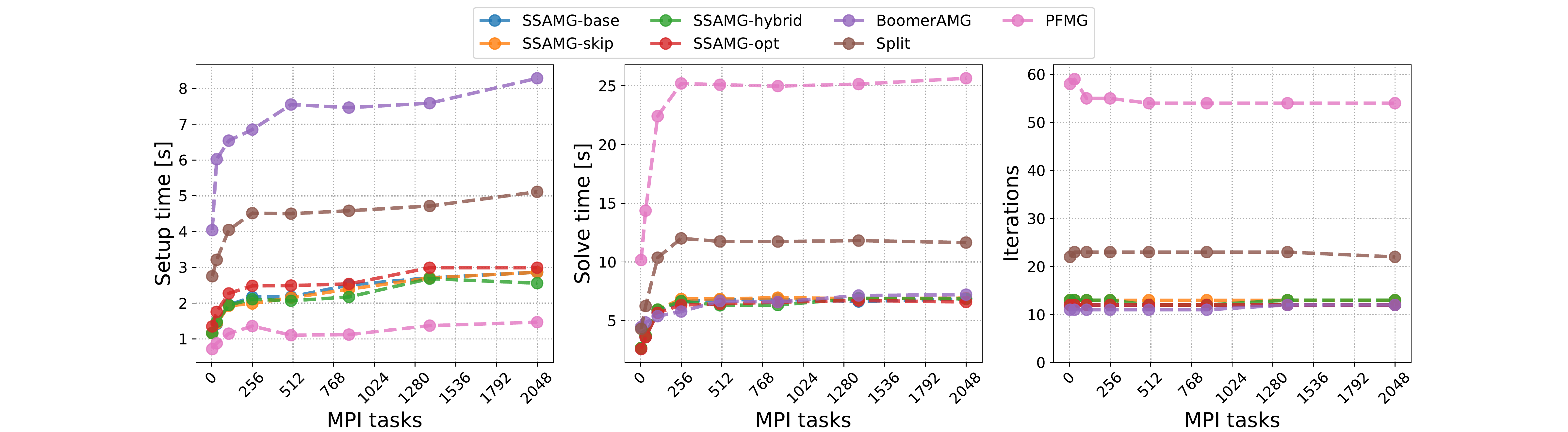}
\caption{Weak scalability results for scenario C of test case 2. Three metrics are shown
  in the figure, i.e., setup phase times in seconds (left); solve phase times in seconds
  (middle), and number of iterations (right). All curves are plotted with respect to the
  number of MPI tasks, $N_{\text{procs}}$, which varies from $4$ ($p = 1$) up to $2048$
  ($p = 8$).}
\label{fig:t1C_ref32x32x32}
\end{figure}

We conclude this section by analyzing the results, given in
Figure~\ref{fig:t1C_ref32x32x32}, for the last anisotropy scenario C. Since there is a
mixed anisotropy configuration in this case as in scenario B, PFMG does not show a
satisfactory convergence behavior. On the other hand, the SSAMG variants show good
numerical and computational scalabilities, and, particularly, \SSAMGOpt{} shows similar
speedups compared to the BoomerAMG variants as discussed in the previous scenarios. When
considering all three scenarios discussed in this section, we note that SSAMG shows good
robustness with changes in anisotropy, and this an important advantage over PFMG.

\subsection{Test case 3 - three-points intersection}
\label{subsec:miller}

\begin{figure}[!hbtp]
\centering
\includegraphics[scale=0.5]{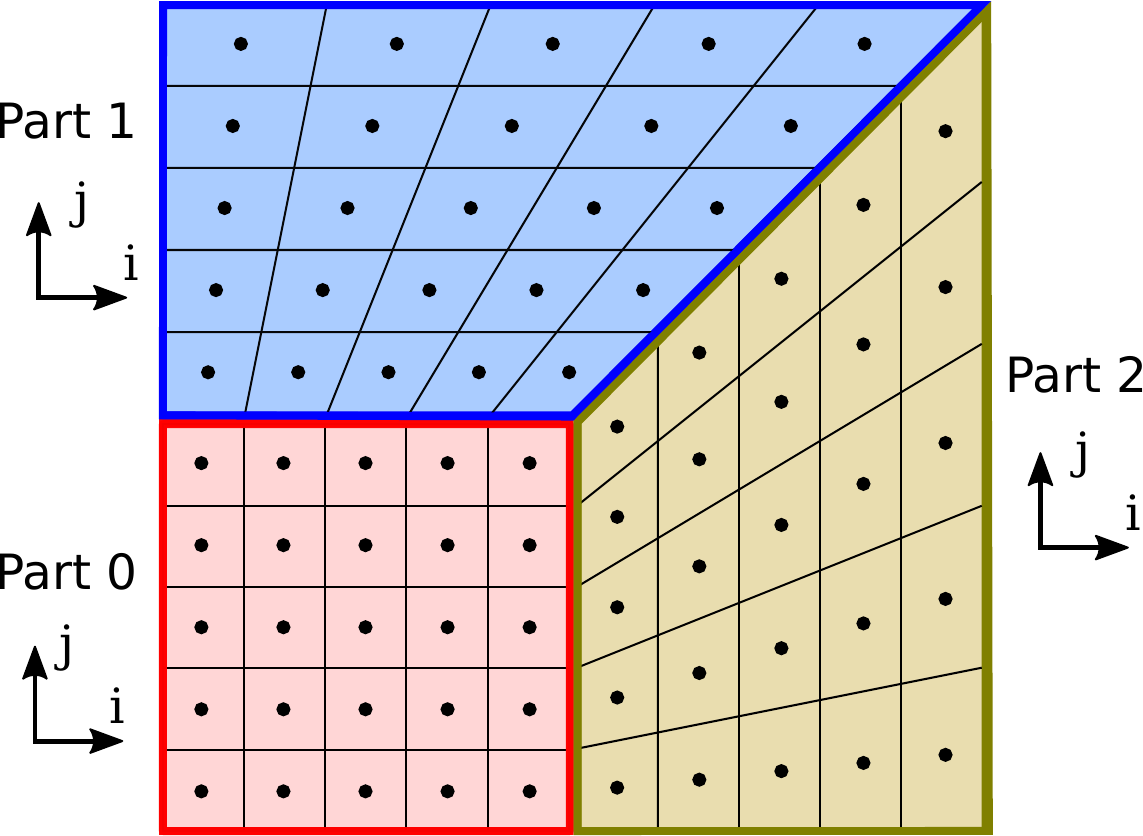}
\caption{$ij$-plane view of the base geometry for test case \ref{subsec:miller}. Equally
  refined instances of this problem in all directions are used for obtaining the results.}
\label{fig:threePoint}
\end{figure}

In this test case, we consider a grid composed topologically of three semi-structured
cubic parts that share a common intersection edge in the $k$-direction
(Figure~\ref{fig:threePoint}). Stencil coefficients are isotropic, but this test case is
globally non-Cartesian. In particular, the coordinate system is different on either side
of the boundary between parts~1 and~2. For example, an east stencil coefficient coupling
Part~1 to Part~2 is symmetric to a north coefficient coupling Part~2 to Part~1.

\begin{figure}[!hbtp]
\centering
\includegraphics[trim={4.5cm 0 4.85cm 0}, clip, width=\textwidth]{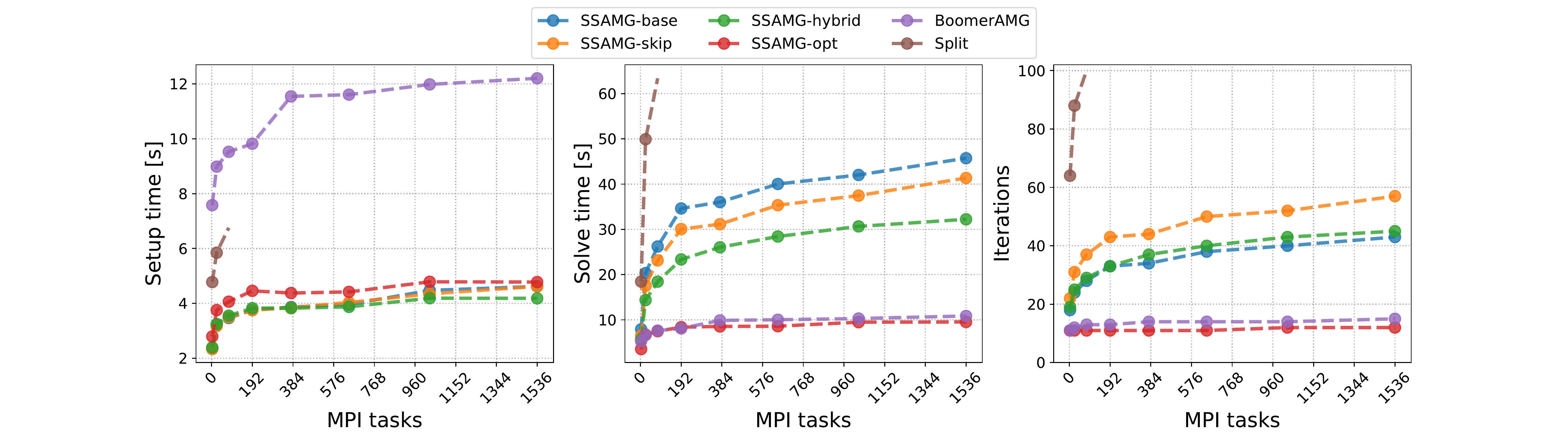}
\caption{Weak scalability results for test case 3. Three metrics are shown in the figure,
  i.e., setup phase times in seconds (left); solve phase times in seconds (middle), and
  number of iterations (right). All curves are plotted with respect to the number of MPI
  tasks, $N_{\text{procs}}$, which varies from $3$ ($p = 1$) up to $1536$ ($p = 8$).}
\label{fig:t2_ref32x32x32}
\end{figure}

For the numerical experiments of this section, we use $m = 160$, which gives a local
problem size per part of $4,096,000$ DOFs, and a global problem size of $12,288,000$ DOFs,
when $p = 1$, i.e., three parts and MPI tasks. Figure~\ref{fig:t2_ref32x32x32} reports
weak scalability results for the current test case. As noted in section
\ref{subsec:cubes}, it is not possible to recast this problem into a single part;
thus, we cannot show results for PFMG here.

Examining the iteration counts reported in Figure~\ref{fig:t2_ref32x32x32}, we see that
\SSAMGOpt{} is the fastest converging option with the number of iterations ranging from $11$, for
$p = 1$ ($3$ MPI tasks), to $14$, for $p = 8$ ($1536$ MPI tasks). This is the best
numerical scalability among the other methods, including BoomerAMG. On the
other hand, the remaining SSAMG variants do not show such good scalability as in
the previous test cases. Once again, this is related to how SSAMG computes interpolation
weights of nodes close to part boundaries. In this context, we plan to investigate further
how to improve SSAMG's interpolation such that the non-hybrid SSAMG variants can have
similar numerical scalability to \SSAMGOpt{}. As in the previous test cases, the
\SStructSplitOpt{} method is the least performing method and does not converge within
100 iterations for $p \geq 3$ ($N_p \geq 81$).

Regarding solve times, \SSAMGOpt{} is the fastest method since it needs the minimum amount
of iterations to reach convergence. Compared to \BoomerAMGOpt{}, its speedup is
$1.3$x for $p = 1$ and $1.1$x for $p = 8$. \SSAMGSkip{} shows solution times smaller than
\SSAMGBase{}, and, here, the ``skip'' option is beneficial to performance. Lastly,
looking at setup times, all SSAMG variants show very similar timings and the optimal variant is up to
$2.9$x faster than \BoomerAMGOpt{}, proving once again the benefits of exploiting problem
structure.

\subsection{Test case 4 - structured adaptive mesh refinement (SAMR)}
\label{subsec:samr}

\begin{figure}[!hbtp]
\centering
\includegraphics[scale=0.5]{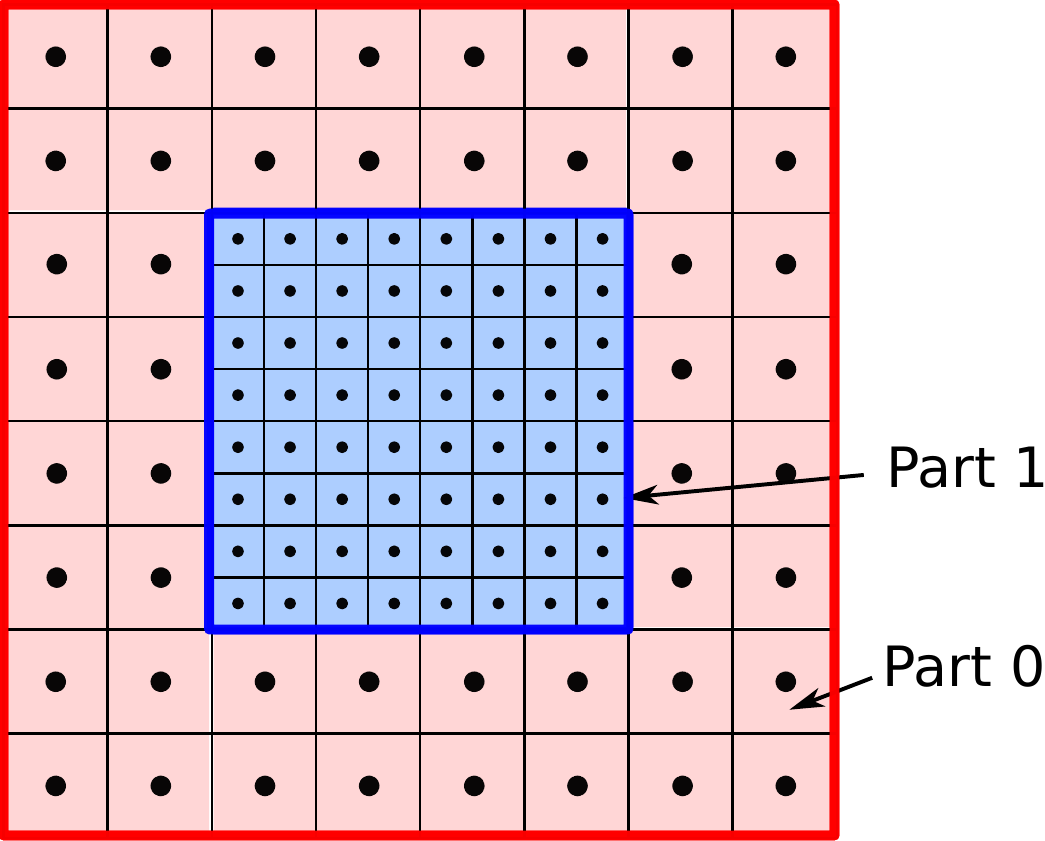}
\caption{$XY$-plane cut of the three-dimensional semi-structured grid used in test case
  \ref{subsec:samr} when $m = 8$. The semi-structured parts represent two levels of
  refinement and contain the same number of cells.}
\label{fig:samr}
\end{figure}

In the last problem, we consider a three-dimensional SAMR grid consisting of one level of grid
refinement, and thus composed of two semi-structured parts (Figure~\ref{fig:samr}). The
first one, in red, refers to the outer coarse grid, while the second, in blue, refers to
the refined patch (by a factor of two) located in the center of the grid. Each part has
the same number of cells. To construct the linear system matrix for this problem,
we treat coarse grid points living inside of the refined part as ghost unknowns, i.e., the diagonal
stencil entry for these points is set to one and the remaining off-diagonal stencil
entries are set to zero. Inter-part couplings at fine-coarse interfaces are stored in the
unstructured matrix ($\mat{U}$), and the value for the coefficients connecting fine grid
cells with its neighboring coarse grid cells (and vice-versa) is set to $2/3$. This value
was determined by composing a piecewise constant interpolation formula with a finite volume
discretization rule. We refer the reader to the SAMR section of \hypre{}'s documentation
\cite{Hypre} for more details.

\begin{figure}[!hbtp]
\centering
\includegraphics[trim={4.8cm 0 4.85cm 0}, clip, width=\textwidth]{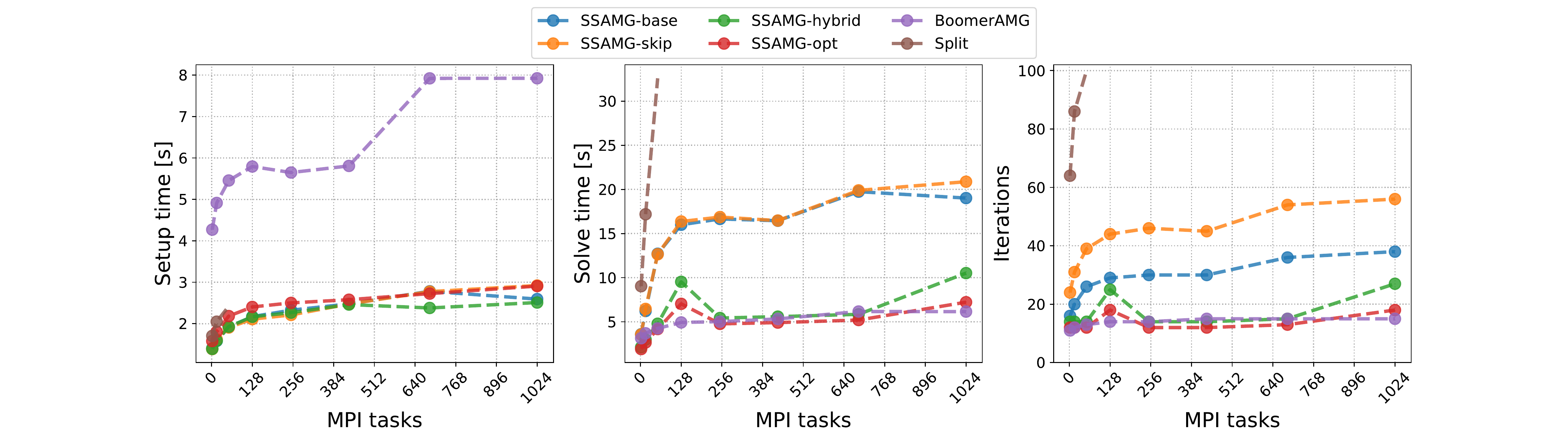}
\caption{Weak scalability results for test case 4. Three metrics are shown in the figure,
  i.e., setup phase times in seconds (left); solve phase times in seconds (middle), and
  number of iterations (right). All curves are plotted with respect to the number of MPI
  tasks, $N_{\text{procs}}$, which varies from $2$ ($p = 1$) up to $1024$ ($p = 8$).}
\label{fig:t3_ref32x32x32}
\end{figure}

The numerical experiments performed in this section used $m = 128$, leading to a local
problem size per part of $2,097,152$ DOFs, and a global problem size of $4,194,304$ DOFs,
for $p = 1$ ($N_{\text{procs}} = 2$). Figure~\ref{fig:t3_ref32x32x32} shows weak
scalability results for this test case. This problem is not suitable for PFMG, thus we do
not show results for PFMG.

As in the previous test cases, \SStructSplitOpt{} does not reach convergence within 100
iterations when $p \geq 3$. Then, \SSAMGSkip{} is the second least convergent option
followed by \SSAMGBase{}. The best option is again \SSAMGOpt{} with the number of
iterations ranging from $15$ ($p = 1$) to $20$ ($p = 8$). Furthermore, its iteration
counts are practically constant for the several parallel runs, except for slight jumps
located at $p = 4$ ($N_{\text{procs}} = 128$) and $p = 8$ ($N_{\text{procs}} = 1024$),
which are present in \SSAMGHybrid{} as well.

As noted before, solve times reflect the methods' convergence performance. In particular,
the cheaper iterations of \SSAMGSkip{} are not able to offset the higher number of
iterations for convergence over \SSAMGBase{}. That explains why these two methods show
very similar solve times.  \SStructSplit{} is the least performing option due to its lack
of robustness. \SSAMGOpt{} and \BoomerAMGOpt{} have similar performance, with \SSAMGOpt{}
slightly better for various cases, but \BoomerAMGOpt{} showing more consistent performance
here.

Setup times of the SSAMG variants are very similar. Listing them in decreasing order of
times, \SSAMGBase{} and \SSAMGSkip{} show nearly the same values, followed by \SSAMGOpt{},
and \SSAMGHybrid{} is the fastest option. The results for \SStructSplitOpt{} are better
here than in the previous test cases, and this is due to the small number of
semi-structured parts involved in this SAMR problem. Still, SSAMG leads to the fastest
options. The slowest method for setup is again \BoomerAMGOpt{} and \SSAMGOpt{} speedups
with respect to it are $4$x for $p = 1$ and $2.7$x for $p = 8$.

\section{Conclusions}
\label{sec:conclusions}

In this paper, we presented a novel algebraic multigrid method, built on the
semi-structured interface in \hypre{}, capable of exploiting knowledge about the problem's
structure and having the potential of being faster than an unstructured algebraic
multigrid method such as BoomerAMG on CPUs and accelerators. Moreover, SSAMG features a
multigrid hierarchy with controlled stencil sizes and significantly improved setup times.

We developed a distributed parallel implementation of SSAMG for CPU architectures in
\hypre{}. Furthermore, we tested its performance, when used as a preconditioner to PCG,
for a set of semi-structured problems featuring distinct characteristics in terms of grid,
stencil coefficients, and anisotropy. SSAMG proves to be numerically scalable for problems
having up to a few billion degrees of freedom and its current implementation achieves
setup phase speedups up to a factor of four and solve phase speedups up to $1.4$x with
respect to BoomerAMG.

For future work, we plan to improve different aspects of SSAMG and its implementation. We
will further investigate SSAMG convergence for more complex problems than have been
considered so far. We want to explore adding an unstructured component to the prolongation
matrix to improve interpolation across part boundaries and evaluate how this benefits
convergence and time to solution. We also plan to add a non-Galerkin option for computing
coarse operators targeting isotropic problems since this approach applied in PFMG has
shown excellent runtime improvements on both CPU and GPU. Finally, we will develop a GPU
implementation for SSAMG.

\section*{Acknowledgments}
This material is based upon work supported by the U.S. Department of Energy, Office of
Science, Office of Advanced Scientific Computing Research, Scientific Discovery through
Advanced Computing (SciDAC) program.

\bibliographystyle{siamplain}
\bibliography{references}

\end{document}